\documentclass{llncs}
\usepackage{amssymb,amsmath,latexsym,epsfig,euscript} 

\newlength{\jmr}
\newlength{\hwl}
\newlength{\khov}
\newlength{\bernd}
\newtheorem{amd}{Arithmetic Multivariate Descartes' Rule (Special Case)} 

\newtheorem{smale}{Smale's $\pmb{\tau}$ Theorem} 

\newtheorem{lenstra}{Lenstra's Theorem}

\newtheorem{thm}{Theorem}
\newtheorem{prop}{Proposition}

\newtheorem{dfn}{Definition}
\newtheorem{cor}{Corollary}

\newtheorem{rem}{Remark}	
 
\newtheorem{ex}{Example}

\newcommand{\np}{{\mathbf{NP}}}

\newcommand{\bpp}{{\mathbf{BPP}}}

\newcommand{\pp}{\mathbf{P}}

\newcommand{\eps}{\varepsilon}

\newcommand{\cO}{O}

\newcommand{\thth}{{\underline{\mathrm{th}}}}

\newcommand{\nd}{{\underline{\mathrm{nd}}}}

\newcommand{\Q}{\mathbb{Q}}
\newcommand{\R}{\mathbb{R}}
\newcommand{\C}{\mathbb{C}}
\newcommand{\N}{\mathbb{N}}
\newcommand{\Z}{\mathbb{Z}}
 
\newcommand{\cp}{\mathfrak{p}}

\newcommand{\Cn}{\C^n}

\newcommand{\Csn}{{(\C^*)}^n}
\renewcommand{\qed}{$\blacksquare$}
\newcommand{\dia}{$\diamond$}

\newcommand{\cL}{{\mathcal{L}}}

\newcommand{\bO}{\mathbf{O}}

\begin{document}

\title{Additive Complexity and Roots of Polynomials Over Number 
Fields and $\cp$-adic Fields} 


\author{J.\ Maurice Rojas \thanks{
This research was partially supported by a grant from the 
Texas A\&M College of Science.}} 

\institute{Department of Mathematics\\ Texas A\&M University\\ TAMU 3368\\
College Station, Texas 77843-3368\\ USA\\ e-mail: {\tt rojas@math.tamu.edu}\\ 
 Web Page: {\tt http://www.math.tamu.edu/\~{}rojas}} 

\date{\today} 

\maketitle

\begin{abstract} 
Consider any nonzero univariate polynomial with rational coefficients,  
presented as an elementary algebraic expression (using only integer 
exponents).  Letting $\sigma(f)$ denotes the additive complexity of $f$, we 
show that the number of rational roots of $f$ is no more than\\ 
\mbox{}\hfill$15+\sigma(f)^2(24.01)^{\sigma(f)} \sigma(f)!$.\hfill\mbox{}\\ 
This provides a sharper 
arithmetic analogue of earlier results of Dima Grigoriev and Jean-Jacques 
Risler, which gave a bound of $C^{\sigma(f)^2}$ for the number of real roots 
of $f$, for $\sigma(f)$ sufficiently large and some constant $C$ with 
$1\!<\!C\!<\!32$. 
We extend our new bound to arbitrary finite extensions of the ordinary {\bf or} 
$p$-adic rationals, roots of bounded degree over a number field, and 
geometrically isolated roots of multivariate polynomial systems. We  
thus extend earlier bounds of Hendrik W.\ Lenstra, Jr.\ and 
the author to encodings more efficient than monomial expansions. 
We also mention a connection to complexity theory and note that our  
bounds hold for a broader class of fields. 
\end{abstract} 

\section{Introduction} 
\label{sec:intro} 
This paper presents another step in the author's program \cite{amd} of 
establishing an effective arithmetic analogue of fewnomial theory. 
(See \cite{few} for the original exposition of 
fewnomial theory, which until now has always used the real or complex numbers 
for the underlying field.) Here, 
we show that the number of {\bf geometrically isolated} roots (cf.\ section 
\ref{sec:multi}) of a polynomial 
system over any fixed $\cp$-adic field (and thereby any fixed number field) 
can be bounded from above by a quantity depending solely on the additive 
complexity of the input equations. 

So let us first clarify the univariate case of {\bf additive complexity}: 
If $\pmb{\cL}$ is any field, we say that $f\!\in\!\cL[x]$ has 
{\bf additive complexity $\pmb{\leq s}$ (over $\pmb{\cL}$)} iff there exist 
constants $c_1,d_1,\ldots,c_s$, $d_s,c_{s+1}\!\in\cL$ and arrays of 
nonnegative integers \scalebox{.92}[1]{$[m_{i,j}]$ and $[m'_{i,j}]$ with 
$f(x)\!=\!c_{s+1}\prod\limits^s_{i=0} X^{m_{i,s+1}}_i$,
where $X_0\!=\!x$, $X_1 = c_1X^{m_{0,1}}_0+
d_1X^{m'_{0,1}}_0$,}\\ and $X_j = c_j \left(\prod\limits^{j-1}_{i=0} 
X^{m_{i,j}}_i\right) + d_j \left(\prod\limits^{j-1}_{i=0} X^{m'_{i,j}}_i 
\right)$ for all $j\!\in\!\{2,\ldots,s\}$. 
We then define the {\bf additive complexity (over $\pmb{\cL}$) of 
$\pmb{f}$}, $\pmb{\sigma_\cL(f)}$, 
to be the least $s$ in such a presentation of $f$ as an algebraic 
expression. Note in particular that additions or subtractions in repeated 
sub-expressions are thus not counted, e.g.,\\ 
$9(x-7)^{99}(2x+1)^{43}-11(x-7)^{999}(2x+1)^3$ has additive complexity 
$\leq\!3$.

It has been known since the work of Allan Borodin and Stephen A.\ Cook 
around 1974 \cite{bocook} that there is a deep connection between 
additive complexity over the real numbers $\R$
and the number of real roots of a nonzero polynomial in $\R[x]$. 
For example, they showed that there is a real constant $K$ such that the 
number of real roots of $f$ is no more than 
$2^{2^{\cdot^{\cdot^{\cdot{^{2^{K\sigma_{_\R}(f)}}}}}}}$, 
where the number of exponentiations is $\sigma_{_\R}(f)-1$ \cite{bocook}. 
Jean-Jacques Risler, using Khovanski's famous Theorem on 
Real Fewnomials \cite{kho,few}, then improved this bound to 
$(\sigma_{_\R}(f)+2)^{3\sigma_{_\R}(f)+1}2^{\left.\left(9\sigma_{_\R}(f)^2 
+5\sigma_{_\R}(f)+2\right)\right/2}$ \cite[pg.\ 181, line 6]{risler}. (Dima 
Grigoriev derived a similar bound earlier \cite{grigo} and both 
results easily imply a simplified bound of $C^{\sigma_{_\R}(f)^2}$ 
for the number of real roots of $f$, for $\sigma_{_\R}(f)$ sufficiently 
large and some constant $C$ with $1\!<\!C\!<\!32$.)  

Here, based on a recent near-optimal {\bf arithmetic} analogue 
of Khovanski's Theorem on Real Fewnomials found by the author (cf.\ 
section \ref{sec:multi} below), we 
give arithmetic analogues of these additive complexity bounds. 
Our first main result can be stated as follows: 
\begin{thm} 
\label{thm:main} 
Let $p$ be any rational prime and let $\log_p(\cdot)$ denote the 
base $p$ logarithm function. Also let $c\!:=\!\frac{e}{e-1}\!\leq\!1.582$, 
let $\cL$ be any degree $d$ algebraic extension of $\Q_p$, and let 
$f\!\in\!\cL[x]\!\setminus\!\{0\}$. 
Then $f$ has no more than 
$2^{\cO\left(\sigma_{_\cL}(f)\log \left(p^d\sigma_{_\cL}(f)\right)\right)}$ 
roots in $\cL$. More precisely, $1+\frac{dp(p^d-1)}{p-1}+
\frac{4cdp(p^d-1)^2}{p-1} \left(1+d\log_p\left(\frac{2d}{\log p}\right)
\right)$\\
\mbox{}\hspace{.5cm} $+\frac{1}{3}\sum\limits^{\sigma_{_\cL}(f)}_{j=3} 
j(6c)^j (p^d-1)^j \left(1+d\log_p\left(\frac{d}{\log 
p}\right)\right) \left(1+d\log_p\left(\frac{2d}{\log p}\right)\right)^{j-1}
j!$ is a valid upper bound, and just the first $\sigma_\cL(f)+1$ 
summands suffice if $\sigma_\cL(f)\!\leq\!2$. 
\end{thm} 
\begin{rem} 
Our bounds can be improved further and this is detailed 
in remark \ref{rem:sharpest} of section \ref{sec:proof}. \dia  
\end{rem} 
\begin{rem}
Note that via the obvious embedding $\Q\!\subset\!\Q_2$, 
theorem \ref{thm:main} easily implies a similar statement for $\cL$ a number 
field. A less trivial extension to number fields 
appears in theorem \ref{thm:bound} below. \dia 
\end{rem} 
\begin{ex} 
\label{ex:first} 
Taking $\cL\!=\!\Q_2$, we obtain respective upper bounds of $1$, $3$, $35$, 
$50195$, \scalebox{.9}[1]{and $6471489$ on the number of roots of $f$ in 
$\Q_2$, according as $\sigma_{_{\Q_2}}(f)$ is $0$, $1$, $2$, $3$, or 
$4$.}\footnote{All calculations in this paper were done with the assistance 
of {\tt Maple} and the corresponding 
{\tt Maple} code can be found on the author's web-page.}\\ 
For instance, we see that for any non-negative integers 
$\alpha,\beta,\gamma,\delta,\eps,
\lambda,\mu,\nu$ and constants $c_1,d_1,c_2,d_2,c_3\!\in\!\Q_2$, 
the polynomial 
\[ c_3x^\alpha\left(c_1x^\beta+d_1x^\gamma\right)^\delta
\left[c_2\left(c_1x^\beta+d_1x^\gamma\right)^\eps+
d_2x^\lambda \left(c_1x^\beta+d_1x^\gamma\right)^\mu\right]^\nu \] 
\scalebox{.86}[1]{has no more than $35$ roots in $\Q_2$ (or $\Q$ obviously).  
See remark \ref{rem:sharper} below for improvements of}\\
\scalebox{.87}[1]{some of these bounds. Note in particular 
that we can {\bf not} count with multiplicities using a}\\ 
\scalebox{.88}[1]{function of $\sigma_\cL(f)$ only, since we can make the 
multiplicities arbitrarily high by increasing}\\
\scalebox{.89}[1]{$\alpha$ and/or $\delta$.  Note 
also that for $\sigma_{_\R}(f)\!\in\!\{0,1,2,3,4\}$ Risler's bound on the 
number of real}\\ 
\scalebox{.91}[1]{roots respectively specializes to $4$, $20736$, 
$274877906944$, $5497558138880000000000$, and}\\ 
\scalebox{.91}[1]{$126315281744229461505151771531542528$.} \dia 
\end{ex} 

The importance of bounds on the number of roots in terms of 
additive complexity is two-fold: on the one hand, we obtain 
a new way to bound the number of roots in $\cL$ of any univariate 
polynomial with coefficients in $\cL$. Going the opposite way, we can use 
information about the number of roots in $\cL$ of a given univariate 
polynomial to give a lower bound on the minimal number of additions and 
subtractions necessary to evaluate it. More to the point, a recent theorem 
of Smale establishes a deep connection between the number of integral  
roots of a univariate polynomial, a variant of additive 
complexity, and certain fundamental complexity classes. 

To make this precise, let us consider another formalization 
of algebraic expressions. Rather than allowing arbitrary recursive 
use of integral powers and field operations, let us be more 
conservative and do the following: 
Suppose we have $f\!\in\!\Z[x_1]$ expressed as a sequence of the form 
$(1,x_1,f_2,\ldots,f_N)$, where $f_N\!=\!f(x_1)$, $f_0\!:=\!1$, 
$f_1\!:=\!x_1$, and for all $i\!\geq\!2$ we 
have that $f_i$ is a sum, difference, or product of some pair of elements 
$(f_j,f_k)$ with $j,k\!<\!i$. (Such computational sequences are 
also known as {\bf straight-line programs} or {\bf SLP's}.) 
Let $\pmb{\tau(f)}$ denote the smallest possible value of $N-1$, i.e., 
the smallest length for such a computation of $f$. Clearly, 
$\tau(f)$ also admits a definition in terms of multivariate polynomial 
systems much like that of $\sigma_\cL(f)$. So it is clear that 
$\tau(f)\!\geq\!\sigma_\cL(f)$ for all $f\!\in\!\Z[x_1]$ and 
$\cL\!\supseteq\!\Z$, and that $\sigma_\cL(f)$ is often dramatically smaller 
than $\tau(f)$. 
\begin{smale}
\cite[theorem 3, pg.\ 127]{bcss} 
Suppose there is an absolute constant $\kappa$ such that 
for all nonzero $f\!\in\!\Z[x_1]$, the 
number of distinct roots of $f$ in $\Z$ is no more than $(\tau(f)+1)^\kappa$. 
Then $\pp_\C\!\neq\!\np_\C$. \qed  
\end{smale} 
In other words, an analogue (regarding complexity 
theory over $\C$) of the famous unsolved $\pp\!\stackrel{?}{=}\!\np$ question 
from computer science (regarding complexity theory over the ring $\Z/2\Z$) 
would be settled. The question of whether $\pp_\C\!\stackrel{?}{=}\!\np_\C$ 
remains open as well but it is known that $\pp_\C\!=\!\np_\C \Longrightarrow 
\np\!\subseteq\!\bpp$. (This observation is due to Steve Smale and 
was first published in \cite{shub}.) The complexity class $\bpp$ is central 
in randomized complexity and cryptology, and the last inclusion (while 
widely disbelieved) is also an open question. The truth of the hypothesis of 
Smale's $\tau$ Theorem, also know as the {\bf $\pmb{\tau}$-conjecture}, is 
yet another open problem, even for $\kappa\!=\!1$. 

Observing that the number of integral roots of 
$f$ is no more than $\deg f$ (by the fundamental theorem of algebra), 
and that $\deg f\!\leq\!2^{\tau(f)}$ (since $\deg f_{i+1}\!\leq\!2\max_{j<i} 
\deg f_j$), we easily obtain the following crude upper bound. 
\begin{prop}
\label{prop:obv} 
\scalebox{.9}[1]{The number of integral roots of 
$f\!\in\!\Z[x_1]\!\setminus\!\{0\}$ is no more than $2^{\tau(f)}$. \qed} 
\end{prop}
As of April 2002, no asymptotically sharper bound in terms of 
$\tau(f)$ appears to be known!\footnote{Using Descartes' Rule of Signs instead 
of the fundamental theorem of algebra does not easily yield a sharper 
bound: the number of monomial terms of $f_i$ grows even faster as a 
function of $\tau(f)$ than $\deg f_i$.} However, taking a $2$-adic approach via 
theorem \ref{thm:main}, we immediately obtain the following improvement.  
\begin{cor}
\label{cor:obv} 
\scalebox{.94}[1]{ The number of integral roots of 
$f\!\in\!\Z[x_1]\!\setminus\!\{0\}$ is 
$2^{^{\cO\left(\sigma_{_{\Q_2}}(f)\log\sigma_{_{\Q_2}}(f)\right)}}$. \qed} 
\end{cor} 
This bound, while apparently not polynomial in $\tau(f)$, 
at least has the advantage that it is frequently much smaller than 
$2^{\tau(f)}$. For instance, our corollary  
tells us that the polynomial from example \ref{ex:first} 
has no more than $35$ integral roots, while the proposition above would 
give us a non-constant upper bound of at least $\alpha$, since this 
example (if not identically zero) has degree $\geq\!\alpha$. 

Whether our $2$-adic approach can be pushed farther to solve 
the $\tau$-conjecture is an intriguing open question. In particular, it 
isn't even known if there is a family of $f$ with 
$2^{^{\Omega\left(\sigma_{_{\Q_2}}(f)\right)}}$ roots in $\Q_2$. 
\begin{rem} 
Curiously, using additive complexity over a different complete field --- $\R$ 
--- can {\bf not} lead to a solution of the $\tau$-conjecture: there are 
examples of $f\!\in\!\Z[x_1]$ 
with $\sigma_{_\R}(f)\!=\!\cO(r)$ and over $2^r$ real (but 
irrational) roots \cite[sec.\ 3, pg.\ 13]{four} (see \cite{bocook} for 
an even bigger lower bound). \dia 
\end{rem} 

Our main results are proved in section \ref{sec:proof}, 
where we in fact prove sharper versions. There we also 
prove a refined number field analogue of theorem 
\ref{thm:main}, which we now state. Recall that if $L$ is a 
subfield of $\C$ and $x\!\in\!\C$ then we say 
\scalebox{.92}[1]{that 
$\pmb{x}$ {\bf is of degree $\pmb{\leq\!\delta}$ over} $\pmb{L}$ iff 
$x$ lies in an algebraic extension of $L$ of degree $\leq\!\delta$.}  
\begin{thm}
\label{thm:bound}
Following the notation of theorem \ref{thm:main}, let 
$\delta\!\in\!\N$ and suppose instead now that $\cL$ is a 
degree $d$ algebraic extension of $\Q$. 
Then the number of roots of $f$ in $\C$ of degree $\leq\!\delta$ over 
$\cL$ is 
$2^{^{\cO\left(\sigma_{_\cL}(f)\left(d\delta+\log \sigma_{_\cL}(f)\right)
\right)}}$. More precisely,\\ 
$1+c(d\delta+10)2^{d\delta+1 }\log_2\left(\frac{d\delta}{\log 2}\right)
+c^2 (d\delta+10)^2 4^{d\delta+2}\log_2\left(\frac{d\delta}
{\log 2}\right)\log_2\left(\frac{2d\delta}
{\log 2}\right)$\\ 
\mbox{}\hspace{.5cm}$+\frac{2}{3} \sum\limits^{\sigma_{_\cL}(f)}_{j=3} 
j(6c)^j 2^{d\delta j} \left(1+2d^2\delta^2\log_2\left(\frac{d^2\delta^2}{\log 
2}\right)\right) \left(1+2d^2\delta^2\log_2\left(\frac{2d^2\delta^2}{\log 2}
\right)\right)^{j-1} j!$\\ is a valid upper bound, and just the 
first $\sigma_\cL(f)+1$ summands suffice if $\sigma_\cL(f)\!\leq\!2$. 
\end{thm} 

\noindent 
This family of bounds can also be sharpened further and this is also 
detailed in remark \ref{rem:sharpest} of section \ref{sec:proof}. 

In summary, theorems \ref{thm:main} and \ref{thm:bound} are  
the first bounds on the number of roots in a local field or 
number field which make explicit use of additive complexity. 
In particular, our results thus extend an earlier result of 
Lenstra on polynomials with few monomial terms to the setting 
of an even sharper input encoding. Recall that for any field 
$L$ we let $L^*\!:=\!L\setminus\{0\}$. 
\begin{lenstra}
\cite[prop.\ 7.2 and prop.\ 8.1]{lenstra2} 
Following the notation of theorems \ref{thm:main} and \ref{thm:bound}, 
suppose now that $\cL$ is a degree $d$ extension of $\Q_p$ (the {\bf local} 
case) {\bf or} $\Q$ (the {\bf global} case), and that 
$f$ has exactly $m$ monomial terms. Then $f$ has no more than 
$c(q_{_\cL}-1)(m-1)^2\left(1+e_{_\cL}
\log_p\left(\frac{e_{_\cL}(m-1)}{\log p}\right)\right)$ roots in $\cL^*$ in the 
local case (counting multiplicities), where $e_{_\cL}$ and $q_{_\cL}$ 
respectively denote the ramification index and residue field cardinality of 
$\cL$.  Furthermore, $f$ has no more than $c(m-1)^2(d\delta+10)\cdot 
2^{d\delta+1} \log_2\left(\frac{d\delta(m-1)}{\log 2}\right)$
roots in $\C^*$ of degree $\leq\!\delta$ over $\cL$ in the global case 
(counting multiplicities). 
\end{lenstra} 
\begin{rem}
Recall that $q_{_\cL}$ is always an integer power of $p$ and 
$e_{_\cL}\log_p q_{_\cL}\!=\!d$. \dia 
\end{rem} 
\begin{ex} 
Considering the polynomial from example \ref{ex:first} once again, note that 
Lenstra's Theorem can not even give a constant upper bound for the number of 
roots in $\Q^*_2$, since the number of monomial terms depends on $\lambda$ 
(among other parameters). On the other hand, in the absence of an  
expression for $f$ more compact than a sum of $m$ monomial terms, Lenstra's 
bound is quite practical. \dia 
\end{ex} 
\begin{rem} 
\label{rem:sharper}
Hendrik W.\ Lenstra has observed that $B(\cL,2,1)$ is in fact the 
number of roots of unity in $\cL$, which is in turn bounded above by 
$\frac{e_{_\cL} p(q_{_\cL}-1)}{p-1}$ \cite{lenstra2}. He has also computed 
$B(\Q_2,3,1)\!=\!6$ (giving $3x^{10}_1+x^2_1-4$ as a trinomial which realizes 
the maximum possible number of nonzero roots in $\Q_2$)  
\cite[prop.\ 9.2]{lenstra2}. As a consequence (following easily from our 
proof of theorem \ref{thm:main}), the first three 
summands of our main formula from theorem \ref{thm:main} 
can be replaced by $1+\frac{e_{_\cL} p(q_{_\cL}-1)}{p-1}+\frac{e_{_\cL} 
p(q_{_\cL}-1)B(\cL,3,1)}{p-1}$, and our bounds from example \ref{ex:first} 
can be improved to $3$ and $15$ in the respective cases 
$\sigma_{_{\Q_2}}(f)\!=\!1$ and 
$\sigma_{_{\Q_2}}(f)\!=\!2$. (This is how we derived the bound cited in the 
abstract.) \dia 
\end{rem} 

As mentioned earlier, our main results follow easily from the 
author's recent arithmetic multivariate analogue of Descartes' Rule 
\cite{amd}. In fact, Arithmetic Multivariate Descates' Rule even allows us to 
derive multivariate extensions of theorems \ref{thm:main} and \ref{thm:bound} 
which we state below. So let us precede our proofs by a brief discussion 
of this important background result. 

\section{Useful Multivariate Results} 
\label{sec:multi}
Suppose $\pmb{f_1},\ldots,\pmb{f_k}\!\in\!\cL[x^{\pm 1}_1,\ldots,
x^{\pm 1}_n]\setminus\{0\}$, and
$\pmb{m_i}$ is the total number of distinct exponent vectors appearing 
in $f_i$ (assuming all polynomials are written as sums of 
monomials). We call $\pmb{F}\!:=\!(f_1,\ldots,f_k)$ a 
$\pmb{k\times n}$ polynomial system over $\pmb{\cL}$ of {\bf type} 
$\pmb{(m_1,\ldots,m_k)}$, and we call 
a root $\zeta$ of $F$ {\bf geometrically isolated} iff $\zeta$ is 
a zero-dimensional component of the underlying scheme over the algebraic 
closure of $\cL$ defined by $F$. If $\cL$ is a finite extension 
of $\Q_p$ (resp.\ $\Q$) then we say that we are in the {\bf local} 
(resp.\ {\bf global}) case. 
\begin{amd}  
\cite[cor.\ 1 of sec.\ 2 and cor.\ 2 of sec.\ 3]{amd} Let $p$ be any 
(rational) prime and $d,\delta$ positive integers. Suppose 
$\cL$ is any degree $d$ algebraic extension of $\Q_p$ or $\Q$, and 
let $\cL^*\!:=\!\cL\setminus\{0\}$. Also let $m\!:=\!(m_1,\ldots,m_n)
\!\in\!\N^n$, $N\!:=\!(N_1,\ldots,N_n)\!\in\!\N^n$, and $F$ an 
$n\!\times\!n$ polynomial system over $\cL$ of type 
$m$ such that the number of variables occuring in $f_i$ is 
exactly $N_i$. Define $\pmb{B(\cL,m,N)}$ to be 
the maximum number of isolated roots in $(\cL^*)^n$ of such an 
$F$ in the local case, counting multiplicities.\footnote{ The 
multiplicity of any isolated root here, which we take in the sense of 
intersection theory for a scheme over the algebraic closure of $\cL$ 
\cite{fulton}, turns out to always be a positive integer when 
$k\!=\!n$ (see, e.g., \cite{smirnov,jpaa}).} 
Then \scalebox{.9}[1]{$\pmb{B(\cL,m,N)\!\leq\!c^n
q^n_\cL \prod\limits^n_{i=1} \left\{m_i (m_i-1)N_i\left[1+e_{_\cL}\log_p
\left(\frac{e_{_\cL}(m_i-1)}{\log p} \right)\right] \right\}}$}, 
where $c\!:=\!\frac{e}{e-1}\!\leq\!1.582$, and $e_{_\cL}$ and $q_{_\cL}$ are 
respectively the ramfication index and residue field cardinality 
of $\cL$. 

Furthermore, moving to 
the global case, let us say a root $x\!\in\!\C^n$ of $F$ is of {\bf degree 
$\pmb{\leq\!\delta}$ over $\cL$} iff every coordinate of $x$  
is of degree $\leq\!\delta$ over $\cL$, and let us define 
$\pmb{A(\cL,\delta,m,N)}$ to be
the maximum number of isolated roots of such an $F$ in $\Csn$ of degree 
$\leq\!\delta$ over $\cL$, counting multiplicities.$^3$ Then 
\scalebox{.98}[1]{$\displaystyle 
\pmb{A(\cL,\delta,m,N)\!\leq\!2c^n2^{d\delta n} \prod\limits^n_{i=1}
\left\{m_i(m_i-1)N_i \left[1+2d^2\delta^2\log_2\left(\frac{d^2\delta^2(m_i-1)}
{\log 2}\right) \right]\right\}}$}. \qed 
\end{amd} 

Various other improvements of these bounds are detailed in \cite{amd}.  
However, let us at least point out that our bound above is nearly optimal: 
For {\bf fixed} $\cL$, $\log B(\cL,(\mu,\ldots,\mu),(n,\ldots,n))$ and 
$\log A(\cL,(\mu,\ldots,\mu),(n,\ldots,n))$ are $\Theta(n\log \mu)$, 
where the implied constant depends on $\cL$ (and $d$ and $\delta$)  
\cite[example 2]{amd}. 

Via our definition of additive complexity we will reduce 
the proofs of our main results to an application of Arithmetic Multivariate 
Descartes' Rule. In particular, it appears that any further improvement to our 
main results will have to come from a different technique. 
For now, we have the following 
generalization of theorems \ref{thm:main} and \ref{thm:bound}. 
\begin{dfn} 
Following the notation above, given any $k\times n$ polynomial 
system $F\!=\!(f_1,\ldots,f_k)$ over $\cL$, 
let us define its {\bf additive complexity over $\pmb{\cL}$}, 
$\pmb{\sigma_\cL(F)}$, to be the 
smallest $s$ such that $F(x_1,\ldots,x_n)$ can be written as\\
$\left(c^{(1)}_{n+s+1}\prod\limits^{n+s}_{i=1} 
X^{m^{(1)}_{i,n+s+1}}_i, \ldots, c^{(k)}_{n+s+1}\prod\limits^{n+s}_{i=1}
X^{m^{(k)}_{i,n+s+1}}_i\right)$, 
where $X_j\!:=\!x_j$ for all\\ $j\!\in\!\{1,\ldots,n\}$,  
$X_j = c_j \left(\prod\limits^{j-1}_{i=1} X^{m_{i,j}}_i\right) + d_j 
\left(\prod \limits^{j-1}_{i=1} X^{m'_{i,j}}_i \right)$ for all\\ 
$j\!\in\!\{n+1,\ldots,n+s\}$, $c_1,d_1,\ldots,c_{n+s},d_{n+s},
c^{(1)}_{n+s+1},\ldots,c^{(k)}_{n+s+1}\!\in\!\cL$,\\ 
and $[m_{i,j}]$, $[m'_{i,j}]$, and $[m^{(\ell)}_{i,j}]$ are arrays of 
positive integers. \dia 
\end{dfn} 
\begin{thm} 
\label{thm:big} 
Following the notation above, $F$ has no more than\\ 
$1+B(\cL,2,1)+(1+B(\cL,2,1)B(\cL,3,1))$\\
\mbox{}\hspace{.5cm}\scalebox{.92}[1]{$+\sum\limits^{\sigma_{_\cL}(F)}_{\ell=3} 
\begin{pmatrix} n+\ell -1\\ n-1 \end{pmatrix}
B(\cL,(\underset{n}{\underbrace{2,\ldots,2}},\underset{\ell-n}
{\underbrace{3,\ldots,3}}),(n+1,n+2\ldots,n+\ell-1,n+\ell-1))$}\\ 
\scalebox{.94}[1]{geometrically isolated roots in $\cL^n$, 
or $1+A(\cL,\delta,2,1)+(1+A(\cL,\delta,2,1)A(\cL,\delta,3,1))$}\\
\mbox{}\hspace{.5cm}\scalebox{.89}[1]{$+\sum\limits^{\sigma_{_\cL}(F)}_{\ell=3} 
\begin{pmatrix} n+\ell -1\\ n-1 \end{pmatrix} 
A(\cL,\delta,(\underset{n}
{\underbrace{2,\ldots,2}},\underset{\ell-n}
{\underbrace{3,\ldots,3}}),(n+1,n+2,\ldots,n+\ell-1,n+\ell-1))$}\\ 
geometrically isolated roots in $\Cn$ of 
degree $\leq\!\delta$ over $\cL$, according as we 
are in the local or global case. 
In particular, for each bound, 
the first $\sigma_\cL(F)+1$ summands suffice if $\sigma_\cL(F)\!\leq\!2$. 
\end{thm} 

In closing, let us point out 
a topological anomaly: Over $\R$, one can go even farther and bound the 
number of connected components of the zero set of a multivariate 
polynomial in terms of additive complexity \cite{grigo,risler}. Unfortunately, 
since $\Q_p$ is totally disconnected as a topological space \cite{koblitz}, 
one can not derive any obvious analogous statement in our 
arithmetic setting. This is why we consider only geometrically 
isolated roots in the multivariate case. Nevertheless, it would be quite 
interesting to know if one could bound the number of higher-dimensional 
{\bf irreducible} components defined over $\cL$ in terms of additive 
complexity, when $\cL$ is a $\cp$-adic field. 

\section{Proving Theorems \ref{thm:main}--\ref{thm:big}} 
\label{sec:proof} 

We will give a proof of Theorem \ref{thm:big} which 
simultaneously yields theorems \ref{thm:main} and 
\ref{thm:bound} for free.  

\noindent
{\bf Proof of Theorem \ref{thm:big} (and Theorems \ref{thm:main} and 
\ref{thm:bound}):} First note that by the definition of additive complexity, 
$(x_1,\ldots,x_n)$ is a geometrically isolated root of $F \Longrightarrow 
(X_1,\ldots,X_{n+s})$ is a geometrically isolated root 
of the polynomial system $G\!=\!\bO$, where the corresponding equations 
are exactly  
{\small  
\[ c^{(1)}_{n+s+1}\prod\limits^{n+s}_{i=1} X^{m^{(1)}_{i,n+s+1}}_i=0 \ \ \ ,  
\ldots, \ \ \ c^{(k)}_{n+s+1}\prod\limits^{n+s}_{i=1}
X^{m^{(k)}_{i,n+s+1}}_i=0,\]
\[X_{n+1} = c_{n+1} \left(\prod\limits^n_{i=1} X^{m_{i,n+1}}_i\right) + d_{n+1}
\left(\prod \limits^n_{i=1} X^{m'_{i,n+1}}_i \right)\] 
\[ \vdots \] 
\[X_{n+s} = c_{n+s} \left(\prod\limits^{n+s-1}_{i=1} X^{m_{i,n+s}}_i\right) + 
d_{n+s} \left(\prod \limits^{n+s-1}_{i=1} X^{m'_{i,n+s}}_i \right),\] 
}

\noindent
where $s\!:=\!\sigma_\cL(F)$, $X_i\!=\!x_i$ for all $i\!\in\!\{1,\ldots,n\}$, 
and the $c_i$, $d_i$, $c^{(j)}_i$, $m_{i,j}$, and $m'_{i,j}$ are 
suitable constants. This follows easily from the fact that 
corresponding quotient rings $\cL[x_1]/\langle f\rangle$ and 
$\cL[X_0,\ldots,X_s]/\langle G \rangle$ are isomorphic, 
thus making $\C_p[x_1]/\langle f\rangle$ and
$\C_p[X_0,\ldots,X_s]/\langle G \rangle$ isomorphic, where 
$\C_p$ denotes the completion of the algebraic closure of $\Q_p$.
In particular, $k\!\leq\!n$ easily implies that 
$F$ has no geometrically isolated roots in $\cL$ at all, 
so we can assume that $k\!\geq\!n$. 

So we now need only count the geometrically isolated roots of 
$G$ in $\cL^{n+s}$ (or the geometrically isolated roots of $F$ in 
$\C^{n+s}$ of degree $\leq\!\delta$ over $\cL$) precisely enough to 
conclude. Toward this end, note that 
the first $n$ equations of \mbox{$G\!=\!\bO$} imply that at 
least $n$ distinct $X_i$ must be $0$, for otherwise 
$(X_1,\ldots,X_{n+s})$ would not be an isolated root. 
Note also that if we have exactly $n$ of the variables 
$X_1,\ldots,X_{n+\ell}$ equal to $0$, then the first $n+\ell$ 
equations of $G$ completely determine $(X_1,\ldots,X_{n+\ell})$. 
Furthermore, by virtue of the last $s-\ell$ equations of $G$, 
the value of $(X_1,\ldots,X_{n+\ell})$ {\bf uniquely} 
determines the value of $(X_{n+\ell+1},\ldots,X_{n+s})$.
So it in fact suffices to find the total number of geometrically isolated 
roots (with all coordinates nonzero) of all systems of the form $G'\!=\!\bO$, 
where the equations of $G'$ are exactly $(0\!=\!0)$ or 
{\small 
\[ \eps_1 X_{n+1} = c_{n+1} \left(\prod\limits^n_{i=1} 
X^{m_{i,n+1}}_i\right) + d_{n+1}
\left(\prod \limits^n_{i=1} X^{m'_{i,n+1}}_i \right)\]
\[ \vdots \]
\[ \eps_\ell X_{n+\ell} = c_{n+\ell} \left(\prod\limits^{n+\ell-1}_{i=1} 
X^{m_{i,n+\ell}}_i\right) + d_{n+s} \left(\prod \limits^{n+\ell-1}_{i=1} 
X^{m'_{i,n+\ell}}_i \right),\] 
}

\noindent
where $\eps_i\!\in\!\{0,1\}$ for all $i$, $X_{n+\ell}\!=\!\eps_\ell\!=\!0$, 
exactly $n-1$ of the variables $X_1,\ldots,X_{n+\ell-1}$ have been set to $0$, 
and $\ell$ ranges over $\{1,\ldots,n\}$. Note in particular that 
the $j^\thth$ equation involves no more than $n+j$ variables for 
all $j\!\in\!\{1,\ldots,\ell-1\}$, and that the $\ell^\thth$ equation 
involves no more than $n+\ell-1$ variables. 

To conclude, we thus see that $G$ has no more than
\[1 \ \ , \ \ 1+B(\cL,2,1) \ \ , \ \ \rho(\cL)\!:=\!1+B(\cL,2,1)+
(r_n+B(\cL,2,1)B(\cL,3,1)) \ \ , \ \ \text{or} \] 
\scalebox{.92}[1]{$\rho(\cL)+\sum\limits^s_{\ell=3} 
\begin{pmatrix} n+\ell -1\\ n-1 \end{pmatrix} 
B(\cL,(\underset{n}
{\underbrace{2,\ldots,2}},\underset{\ell-n}
{\underbrace{3,\ldots,3}}),(n+1,n+2,\ldots,n+\ell-1,n+\ell-1))$}\\ 
geometrically isolated roots in $\cL^{n+s}$ in the 
local case, according as $s$ is $0$, $1$, $2$, or $\geq\!3$, 
where $r_n$ is $0$ or $1$ according as $n\!=\!1$ or $n\!\geq\!2$.  
The corresponding statement for the global case, where we replace $B(\cL,m,N)$ 
by $A(\cL,\delta,m,N)$ throughout and count geometrically isolated roots in 
$\C^{n+s}$ of degree $\leq\!\delta$ over $\cL$ instead, is also clearly true. 
This proves theorem \ref{thm:big}. 

Theorems \ref{thm:main} and \ref{thm:bound} then follow immediately by 
specializing the above formulae to $n\!=\!1$, applying Arithmetic 
Multivariate Descartes' Rule, and performing an elementary calculation. \qed 

\begin{rem} 
\label{rem:sharpest} 
It follows immediately from our proof that we can restate 
theorems \ref{thm:main} and \ref{thm:bound} in sharper 
intrinsic terms. That is, the bounds from our proof above 
can immediately incorporate any new upper bounds for 
the quantities $B(\cL,m,N)$ and $A(\cL,\delta,m,N)$. \dia 
\end{rem} 

\begin{rem} 
Note that the same proof will essentially work verbatim if 
we replace \scalebox{.95}[1]{$\cL$ throughout by {\bf any} field admitting a 
multivariate analogue of Descartes' Rule. \dia}\footnote{In 
particular, via the approach of our proofs, it is possible to improve slightly 
the bounds of \cite{grigo,risler} over $\R$. We leave this as an exercise for 
the interested reader.} 
\end{rem} 

\noindent
{\bf Acknowledgement} The author thanks the two anonymous referees 
for their astute comments. 

\footnotesize
\bibliographystyle{acm}

\end{document}